\pdfoutput=1
\RequirePackage{ifpdf}
\ifpdf 
\documentclass[pdftex]{sigma}
\else
\documentclass{sigma}
\fi

\newtheorem{prop}{Proposition}[section]
\newtheorem{thm}[prop]{Theorem}

\newtheorem{lem}[prop]{Lemma}

\theoremstyle{definition}
\newtheorem{defn}[prop]{Definition}

\newtheorem{rem}[prop]{Remark}

\numberwithin{equation}{section}

\DeclareMathOperator\kernel{ker}
\DeclareMathOperator\sym{sym}

\DeclareMathOperator\Hom{Hom}

\DeclareMathOperator\GL{GL}
\DeclareMathOperator\SL{SL}
\DeclareMathOperator\Diff{Dif\/f}
\DeclareMathOperator\SDiff{SDif\/f}
\DeclareMathOperator\HDiff{HDif\/f}
\DeclareMathOperator\Aff{Af\/f}

\begin{document}
\newcommand{\ip}[1]{\langle#1\rangle}
\newcommand{\abrack}[1]{[\mkern-3mu[#1]\mkern-3mu]}
\newcommand{\vspan}[1]{\left<#1\right>}

\allowdisplaybreaks

\renewcommand{\thefootnote}{$\star$}

\renewcommand{\PaperNumber}{035}

\FirstPageHeading

\ShortArticleName{Selfdual 4-Manifolds, Projective Surfaces, and the Dunajski--West Construction}

\ArticleName{Selfdual 4-Manifolds, Projective Surfaces,\\
and the Dunajski--West Construction\footnote{This paper is a~contribution to the Special Issue on Progress in Twistor
Theory. The full collection is available at
\href{http://www.emis.de/journals/SIGMA/twistors.html}{http://www.emis.de/journals/SIGMA/twistors.html}}}

\Author{David M.J.~CALDERBANK}
\AuthorNameForHeading{D.M.J.~Calderbank}
\Address{Department of Mathematical Sciences, University of Bath, Bath BA2 7AY, UK}
\Email{\href{mailto:D.M.J.Calderbank@bath.ac.uk}{D.M.J.Calderbank@bath.ac.uk}}
\URLaddress{\url{http://people.bath.ac.uk/dmjc20/}}

\ArticleDates{Received January 21, 2014, in f\/inal form March 18, 2014; Published online March 28, 2014}

\Abstract{I present a~construction of real or complex selfdual conformal $4$-manifolds (of signature $(2,2)$ in the real
case) from a~natural gauge f\/ield equation on a~real or complex projective surface, the gauge group being the group of
dif\/feomorphisms of a~real or complex $2$-manifold.
The $4$-manifolds obtained are characterized by the existence of a~foliation by selfdual null surfaces of a~special
kind.
The classif\/ication by Dunajski and West of selfdual conformal $4$-manifolds with a~null conformal vector f\/ield is the
special case in which the gauge group reduces to the group of dif\/feomorphisms commuting with a~vector f\/ield, and I
analyse the presence of compatible scalar-f\/lat K\"ahler, hypercomplex and hyperk\"ahler structures from
a~gauge-theoretic point of view.
In an appendix, I discuss the twistor theory of projective surfaces, which is used in the body of the paper, but is also
of independent interest.}


\Classification{53A30; 32L25; 37K25; 37K65; 53C25; 70S15; 83C20; 83C60}

\renewcommand{\thefootnote}{\arabic{footnote}}
\setcounter{footnote}{0}

\section{Introduction}

In~\cite{DuWe:acs}, M.~Dunajski and S.~West presented a~local classif\/ication of those selfdual $4$-manifolds, either
complex or real with signature $(2,2)$ (known also as kleinian, ultrahyperbolic, neutral or split signature), that admit
a~nontrivial null conformal vector f\/ield.
To do this, they observed that the two natural null plane distributions associated to a~null conformal vector f\/ield are
integrable.
By exploiting one of the two integral foliations, they reduced the geometry to that of a~projective surface.

The two null plane distributions have opposite orientations with respect to the Hodge star operator, one being selfdual,
the other antiselfdual.
Such null planes (or their integral surfaces) are referred to as \textit{$\alpha$-planes} and \textit{$\beta$-planes}
(or surfaces), following the Klein correspondence between ${\mathbb{R}}{\rm P}^3$ or ${\mathbb{C}}{\rm P}^3$
and the real or complex Klein quadrics in ${\rm P}(\wedge^2{\mathbb{R}}^4)$ or
${\rm P}(\wedge^2{\mathbb{C}}^4)$, which are the four dimensional conformally-f\/lat spaces of lines in these
projective spaces.
The Klein quadrics contain two classes of null projective planes called $\alpha$-planes and $\beta$-planes: the
$\alpha$-planes correspond to (the lines through) points in ${\mathbb{R}}{\rm P}^3$ or
${\mathbb{C}}{\rm P}^3$, whereas the $\beta$-planes correspond to (the lines in) planes in
${\mathbb{R}}{\rm P}^3$ or ${\mathbb{C}}{\rm P}^3$.

The selfduality condition on a~conformal $4$-manifold $M$ is the integrability condition for the existence of an
antiselfdual surface tangent to a~given antiselfdual plane at any point.
The space of such surfaces, for suitably convex $M$, is then a~$3$-manifold generalizing ${\mathbb{R}}{\rm P}^3$ or
${\mathbb{C}}{\rm P}^3$, called the \textit{twistor space} $Z$, and so antiselfdual $2$-planes are called
$\alpha$-planes, whereas selfdual $2$-planes are called $\beta$-planes.
Using the identif\/ication $T_xM\cong S_x\otimes S'_x$, where $S$ and $S'$ are the (rank $2$) selfdual and antiselfdual
spinor bundles, $\alpha$-planes at $x$ are $2$-planes of the form $S_x\otimes \vspan{\ell'_x}$ for $\vspan{\ell'_x}\in
{\rm P}(S'_x)$, while the $\beta$-planes instead have the form $\vspan{\ell_x}\otimes S'_x$ for $\vspan{\ell_x}\in
{\rm P}(S_x)$.
(Dunajski and West use the opposite orientation, so that their $\alpha$-surfaces are selfdual.
I have chosen to retain the term $\alpha$-surface for the null surfaces corresponding to points in twistor space, even
though these are antiselfdual for my choice of orientation.)

A point $x\in M$ may be identif\/ied with the set of $\alpha$-surfaces which pass through it, and this def\/ines a~rational
curve $\hat x$ in $Z$, called a~\textit{twistor line}, which has normal bundle $\mathcal{O}_{\hat x}(1)\otimes S_x$ in
$Z$, where $\mathcal{O}_{\hat x}(1)$ is the dual of the tautological line bundle of $\hat x$.
In the complex case, this leads to a~construction of selfdual conformal manifolds as Kodaira moduli spaces of such
twistor lines~\cite{Pen:nlg}.

The existence of an $\alpha$-surface foliation on a~selfdual conformal manifold is therefore not remarkable: such
foliations simply correspond to hypersurfaces in $Z$ transverse to twistor lines.
On the other hand, the existence of a~$\beta$-surface foliation constrains the geometry signif\/icantly.
It is the $\beta$-surface foliation that Dunajski and West exploit: they show that the (locally def\/ined) leaf space of
$\beta$-surfaces is a~projective surface, and use this to completely classify (anti)selfdual
conformal $4$-manifolds with a~null conformal vector f\/ield.

My purpose in this paper is to determine when a~$\beta$-surface foliation induces a~projective structure on the leaf
space, and to show that there is a~generalized Dunajski--West construction of any selfdual conformal $4$-manifold, with
such a~foliation, from a~natural gauge f\/ield equation on the projective surface.
I show how the Dunajski--West classif\/ication f\/its into this construction as a~reduction of gauge group, and analyse the
presence of compatible compatible scalar-f\/lat K\"ahler, hypercomplex and hyperk\"ahler structures from a~gauge-theoretic
point of view.

I make use of the twistor theory of projective surfaces introduced by C.~LeBrun~\cite{LeB:phd} and
N.~Hitchin~\cite{Hit:cme,Hit:gse} (see also~\cite{Dub:scv}).
Since this has independent interest, and has not yet been developed in detail, I discuss some aspects of this twistor
theory, including Penrose transforms, Ward correspondences, and divisors, in an appendix.
(In the projectively-f\/lat case, the Penrose transform has already been used to provide a~twistorial approach to the Funk
and Radon transforms~\cite{BEGM:ftpt}.)

Since this article f\/irst appeared on the arXiv, F.~Nakata has developed a~global approach~\cite{Nak:sdza} relating
$\beta$-surface reduction to the work of C.~LeBrun and L.~Mason~\cite{LeMa:zmcs,LeMa:nghd}.

\section{The twistor correspondence for beta surface foliations}

For complex (i.e., holomorphic) conformal $4$-manifolds (and hence, via complexif\/ication, for real analytic
conformal $4$-manifolds of signature $(2,2)$), there is a~simple twistorial approach.
First recall~\cite{Hit:cme,Hit:gse} that there is twistor-like correspondence between complex projective surfaces~$N$
and complex surfaces~$Y$ containing rational curves with normal bundle $\cong\mathcal{O}(1)$, which I will call
\textit{minitwistor spaces} (and the rational curves will be called \textit{minitwistor lines}).
Given a~geodesically convex complex projective surface~$N$, the space of (unparameterized) geodesics is a~minitwistor
space $Y$, the minitwistor lines corresponding to the (geodesics through) points in~$N$.
Conversely, for any minitwistor space $Y$, the moduli space of its minitwistor lines is, by Kodaira deformation theory,
(locally) a~complex surface, and the paths corresponding to points in~$Y$ are the geodesics of a~projective structure.

{\sloppy A special case of this correspondence is the correspondence between dual projective pla\-nes~${\rm P}({\mathbb{C}}^3)$ and ${\rm P}({\mathbb{C}}^{3*})$, or, more generally, between a~convex open subset~$N$ of $P({\mathbb{C}}^3)$ and the set~$Y$ of points in ${\rm P}({\mathbb{C}}^{3*})$ corresponding to lines that
meet~$N$.}

Now suppose that $Y$ is a~minitwistor space over which there is a~f\/ibre bundle $Z\to Y$ of rank one and degree one: by
this I mean that the pullback of $Z$ to a~minitwistor line has a~section whose normal bundle (which is the vertical
bundle of $Z$ along the section) is $\cong\mathcal{O}(1)$; in particular $Z$ has one-dimensional f\/ibres.
Then, by Kodaira deformation theory, $Z$ (the total space) is a~twistor space: the sections of $Z$ along a~minitwistor
line are rational curves with normal bundle $\cong\mathcal{O}(1)\otimes{\mathbb{C}}^2$.
The twistor spaces $Z$ arising in this way are those that admit a~foliation by curves which are transverse to its
twistor lines and such that a~twistor line $\hat x$ has normal bundle $\cong\mathcal{O}(1)$ inside the union $H_x$ of
the leaves meeting it: $Y$ is then the (local) leaf space of this foliation, the normal bundles in $H_x$, $Z$ and $Y$
being related by the obvious exact sequence
\begin{gather*}
0\to \mathcal{O}(1)\to \mathcal{O}(1)\otimes{\mathbb{C}}^2\to \mathcal{O}(1)\to 0
\end{gather*}
(which is the only extension of $\mathcal{O}(1)$ by $\mathcal{O}(1)$ up to isomorphism).
Such a~foliation determines a~line subbundle of $TZ$ of degree one (i.e., restricting to a~degree one line
bundle on each twistor line).
Let $\mathcal{O}_Z(1)$ be the fourth root $K_Z^{-1/4}$ of the anticanonical bundle of $Z$ (which has degree one).
There is thus a~line subbundle
\begin{gather*}
\hat\ell\colon \ \mathcal{L}\otimes\mathcal{O}_Z(1)\hookrightarrow TZ,
\end{gather*}
where $\mathcal{L}$ is degree zero (i.e., trivial on twistor lines).
Conversely such a~line subbundle integrates to a~foliation of the desired type: it is transverse to generic twistor
lines because the projection of $\hat\ell$ onto the normal bundle of a~twistor line $\hat x$ is a~section of a~bundle
isomorphic to $\Hom(\mathcal{O}(1),\mathcal{O}(1)\otimes{\mathbb{C}}^2)$ along $\hat x$, hence constant.

The f\/irst task is to determine the object corresponding to $\hat\ell$ on the space $M$ of twistor lines in $Z$ to which
$\hat\ell$ is transverse.
This is a~straightforward exercise in computing the Penrose transform.
First restrict $\hat\ell$ to a~twistor line $\hat x$ corresponding to $x\in M$ to obtain a~(nonzero) section of
$\Hom(\mathcal{L}\otimes\mathcal{O}_Z(1)|_{\hat x},\mathcal{O}_{\hat x}(1)\otimes S_x)$.
Let~$L_x$ be the space of sections of $\mathcal{L}\otimes\mathcal{O}_Z(1)|_{\hat x}\otimes\mathcal{O}_{\hat x}(-1)$,
which is a~complex line; then $\hat\ell$ induces an inclusion $\ell_x\colon L_x\hookrightarrow S_x$, hence
a~$\beta$-plane in~$T_x M$.
Varying~$x$ yields a~line subbundle $L$ of $S$ over $M$, and hence a~$\beta$-plane distribution~$\ell$.
I will show that it is in fact a~$\beta$-surface foliation.
However, not every $\beta$-surface foliation arises in this way.
For this suppose f\/irst that~$\mathcal{L}$ is trivial and adopt the convention that~$\mathcal{O}_{\hat
x}(1)=\mathcal{O}_Z(1)|_{\hat x}$ for all~$x$: this amounts to requiring that $S$ is the weight $1/2$ selfdual spinor
bundle, on which the twistor operator $\mathcal{T}$ is conformally invariant, and it turns out that $\ell$ is a~twistor
spinor, i.e., $\mathcal{T}\ell=0$.
More generally, if $\mathcal{L}$ is not trivial it corresponds to a~connection $\nabla$ on $L$ with selfdual curvature
(also known as a~selfdual Maxwell f\/ield), and $\mathcal{T}^\nabla\ell=0$.

Now consider an arbitrary $\beta$-plane distribution $\ell\colon L\hookrightarrow S$.
Then for any connection~$\nabla$ on~$L$, one can def\/ine a~section $ \mathcal{T}^\nabla\ell$ of $T^{*} M\odot S\otimes
L$, where $T^{*} M\odot S:=\kernel c\colon T^{*} M\otimes S \to S'$ and~$c$ is Clif\/ford multiplication, by projecting
$D^\nabla\ell$ onto this kernel, where $D$ is the Levi-Civita connection of any compatible metric, or indeed any Weyl
connection.
(In spinor index notation $(\mathcal{T}^\nabla\ell)_{A'AB}=D^\nabla_{A'(A}\ell_{B)}$.) As is well known
(see~\cite{DuWe:acs,DuWe:acns}), $\ell$ is integrable, i.e., tangent to a~$\beta$-surface foliation, if and
only if $\ip{\mathcal{T}^\nabla\ell, \ell\otimes\ell}=0$ (i.e.,
$(\mathcal{T}^\nabla\ell)_{A'AB}\ell^A\ell^B=0$), a~condition which is independent of~$\nabla$, since
$\mathcal{T}^{\nabla+\gamma}\ell=\mathcal{T}^\nabla\ell+\gamma\odot\ell$ for any $1$-form $\gamma$ (in indices
$\gamma\odot\ell$ is $\gamma_{A'(A}\ell_{B)}$).
However this integrability condition implies that $\mathcal{T}^\nabla\ell$ is of the form $\alpha\odot\ell$ for
a~$1$-form $\alpha$ and hence there is always a~connection $\nabla$ such that $\mathcal{T}^\nabla\ell=0$.
Furthermore $\gamma\odot \ell=0$ implies $\gamma=0$ and so the connection is unique.
I therefore refer to $\nabla$ as the canonical connection of the $\beta$-surface foliation and make the following
def\/inition.

\begin{defn}
A $\beta$-surface foliation $\ell\colon L\hookrightarrow S$ with canonical connection $\nabla$ on $L$ will be called
\textit{selfdual} if $\nabla$ has selfdual curvature.
\end{defn}

Lest any confusion arise, this condition has nothing to do with the fact that $\beta$-planes are selfdual $2$-planes:
selfdual $\beta$-surface foliations are, in a~sense, ultra-selfdual! The terminology follows the paper~\cite{Cal:sde},
where I introduced the notion of a~selfdual conformal submersion, which involves a~similar selfduality property for line
subbundles of~$TM$.

I hope that the following result does not now come as a~surprise.

\begin{thm}
\label{main}
Let $M$ be a~complex selfdual conformal $4$-manifold with twistor space $Z$.
Then there is a~local correspondence between selfdual $\beta$-surface foliations of $M$ and foliations of $Z$ by curves
which are transverse to twistor lines, and such that the twistor lines have normal bundle~$\mathcal{O}(1)$ in the union
of the curves which meet them.

Hence if $M$ is such a~manifold with a~selfdual $\beta$-surface foliation, the local leaf space~$Y$ of the corresponding
curves in~$Z$ is a~minitwistor space.
The local leaf space~$N$ of $\beta$-surfaces is the corresponding projective surface.
\end{thm}

Note that this is a~statement in the holomorphic category, i.e., the $\beta$-surface foliations and foliations
of $Z$ are assumed to be holomorphic.
It is also local to a~neighbourhood of a~point in~$M$ and the corresponding twistor line in~$Z$.
This permits the assumption in the argument below that intersections are generic, and that constructed objects are
nonsingular.

\begin{proof}
Suppose f\/irst that $Z$ has a~foliation by curves as described, let $\hat x$ be a~twistor line corresponding to $x\in M$,
and consider the union $H_x$ of the curves in $Z$ which meet $\hat x$.
This is a~(locally nonsingular) divisor in $Z$, inside which $\hat x$ has normal bundle $\cong \mathcal{O}(1)$.
As shown above, a~tangent vector f\/ield to these curves along $\hat x$ projects onto a~section of $\mathcal{O}_{\hat
x}(1)\otimes\vspan{s}\subset \mathcal{O}_{\hat x}(1)\otimes S_x$ (the normal bundle of $\hat x$ in $Z$) for some nonzero
$s\in S_x$ uniquely determined up to scale, corresponding to a~$\beta$-plane at $x\in M$.

Now by Kodaira deformation theory, the moduli space of deformations of $\hat x$ in $H_x$, which is (locally)
a~submanifold of $M$, has tangent space $H^0(\hat x,\mathcal{O}_{\hat x}(1))\otimes \vspan{s}$, which is the
corresponding $\beta$-plane.
Hence it is a~surface $\Sigma$ in $M$ tangent to the given $\beta$-plane at $x$.
The same applies to any $x_1\in\Sigma$ (since $Z$ is foliated, $H_{x_1}=H_x$), so that $\Sigma$ is a~$\beta$-surface.
Now, varying $x\in M$ yields an integrable $\beta$-plane distribution, hence a~foliation by $\beta$-surfaces.

Let $\hat\ell\colon \mathcal{L}\otimes \mathcal{O}_Z(1)\hookrightarrow TZ$ and $\ell\colon L\hookrightarrow S$ be
tangent to the corresponding foliations, as in the discussion above.
By the Penrose--Ward transform for selfdual Maxwell f\/ields (a special case of the Ward correspondence also known as the
twisted photon construction), $\mathcal{L}$ induces on $L$ a~connection $\nabla$ with selfdual curvature whose parallel
sections along any $\alpha$-surface $\hat z$ are precisely those evaluating (at $x\in \hat z$) to a~f\/ixed element of
$\mathcal{L}_z$ (when viewed as a~section of $\mathcal{L}$ over $\hat x$).
Now restrict $\ell$ to $\hat z$ and trivialize $L$ by parallel sections corresponding to $e\in\mathcal{L}_z$.
Then by construction, for any $x\in\hat z$, $\ell_x\otimes e=\hat\ell_z \mathbin{{\rm mod}} T\hat x$.
It follows by a~standard twistor theory that the contraction of $\mathcal{T}^\nabla\ell$ with the antiselfdual spinor
corresponding to $z$ vanishes along $\hat z$.
Since $\hat z$ was arbitrary $\mathcal{T}^\nabla\ell=0$.
Hence $\nabla$ is the canonical connection of $\ell$, and $\ell$ is selfdual.
(In fact this gives another proof that the $\beta$-plane distribution is integrable, but I f\/ind the Kodaira argument
more conceptually appealing.)

For the converse, suppose $\ell$ is a~selfdual $\beta$-surface foliation with canonical connection $\nabla$ and let
$\Sigma$ be a~surface in the foliation.
The set of $\alpha$-surfaces meeting $\Sigma$ is a~(locally nonsingular) divisor $\hat\Sigma$ in $Z$: indeed such
$\alpha$-surfaces are uniquely determined (holomorphically) by the tangent line to their intersection with $\Sigma$ at
some point, and so form a~$2$-parameter family (because they meet $\Sigma$ in a~curve).
Now $x\in\Sigma$ corresponds to a~twistor line $\hat x$ contained in $\hat\Sigma$ (parameterized by the tangent lines in
$T_x \Sigma$), hence $\hat\Sigma$ contains a~$2$-parameter family of twistor lines.
A vector f\/ield tangent to $\hat\Sigma$ along $\hat x$ clearly projects to a~section of $\mathcal{O}_{\hat
x}(1)\otimes\ell_x$.
It follows that the twistor lines in $\hat\Sigma$ have normal bundle $\mathcal{O}(1)$ and that $\Sigma$ is (locally)
uniquely determined by $\hat\Sigma$ as the moduli space of twistor lines in $\hat\Sigma$.

Now f\/ix an $\alpha$-surface $\hat z$ corresponding to $z\in Z$.
Since $\nabla$ has selfdual curvature, it is f\/lat along $\hat z$, and the space of parallel sections along $\hat z$ is
the f\/ibre $\mathcal{L}_z$ of a~line bundle $\mathcal{L}$ on $Z$ (this is the inverse Penrose--Ward transform).
On the other hand $\mathcal{T}^\nabla\ell$ vanishes along $\hat z$.
Hence on trivializing $\mathcal{L}_z$ by an element $e$ and $L|_{\hat z}$ by the corresponding parallel section,
$\ell_x\otimes e = \hat\ell_z\mathbin{{\rm mod}} T\hat x$ for all $x\in\hat z$, where $\hat\ell_z$ is a~f\/ixed
element of $\Hom(\mathcal{L}\otimes\mathcal{O}_Z(1),TZ)_z$, independent of $x\in\hat z$ and the choice of $e$.
Varying $z$, thus yields a~section of $\Hom(\mathcal{L}\otimes\mathcal{O}_Z(1),TZ)$, def\/ining a~foliation by curves.
By construction, the foliation has the required properties.
It is also easy to see that the two constructions are mutually inverse.

As observed in~\cite{DuWe:acs}, the geometry of the divisors $\hat\Sigma$, corresponding to the $\beta$-surfaces
$\Sigma$ in the foliation is quite easy to understand.
For any $z\in Z$, the $\beta$-surfaces in the foliation which meet~$\hat z$ are determined by their intersections with
$\hat z$, which foliate it, hence they form a~$1$-parameter family~$\Sigma_t$, corresponding to a~$1$-parameter family
of distinct divisors $\hat\Sigma_t$ through $z\in Z$.
The tangent spaces to these divisors all contain the direction def\/ined by $\hat\ell_z$.
Since this is true for all~$z$, it follows that each curve in the foliation of $Z$ is the intersection of
a~$1$-parameter family of $\hat\Sigma$'s, so that it def\/ines $1$-parameter families of $\alpha$ and $\beta$-surfaces
meeting each other pairwise in curves.

Therefore, the correspondence descends to the local leaf spaces $N$ and $Y$ of the foliations: a~point in~$N$ is
a~$\beta$-surface $\Sigma$, corresponding to a~divisor $\hat\Sigma$ in $Z$, which descends to a~minitwistor line in~$Y$;
inversely, a~point in $Y$ def\/ines a~curve in $Z$, hence a~$1$-parameter family of $\beta$-surfaces in $M$, which is
a~path in $N$.
This completes the proof.
\end{proof}

It is easy to check that these constructions can be done compatibly with real structures for the spaces involved (where
the real points in $M$ form a~kleinian signature conformal manifold).
This settles the real analytic case.
The extension to kleinian signature conformal manifolds which are merely smooth, requires further work.
One approach is to use the generalized Kodaira theory of C.~LeBrun and L.~Mason~\cite{LeMa:nghd}.
However, following Dunajski and West, I will instead construct explicitly the correspondence spaces, since this approach
yields explicitly the form of the construction of $M$ from a~gauge f\/ield equation on $N$.

\section{The inverse construction from projective pairs}

So far, the only procedure I have given for constructing a~selfdual conformal manifold $M$ from a~projective surface $N$
is rather indirect: $M$ is a~Kodaira moduli space of curves in a~rank one, degree one f\/ibre bundle $Z\to Y$.
However, Kodaira deformation theory can almost never be carried out explicitly, so a~more direct approach is required.
To motivate a~more explicit inverse construction, and the generalization of the theory to the smooth case, I now look at
the structure of the correspondence space revealed by the Theorem~\ref{main} when $M$ and $Z$ are complex.
In general, the correspondence space is the manifold
\begin{gather*}
F=\{(x,z)\in M\times Z:x\in \hat z\}.
\end{gather*}
The obvious projections $\pi_M$ and $\pi_Z$ make $F$ into a~${\mathbb{C}}{\rm P}^1$ bundle over $M$ and (without
loss, since the construction is local on $M$) a~${\mathbb{C}}^2$ bundle over $Z$.
As a~bundle over $M$, $F=PS'$, the projectivized antiselfdual spinor bundle, and its f\/ibres project to twistor lines in
$Z$.
The f\/ibration $\pi_Z\colon F\to Z$ (whose f\/ibres project to $\alpha$-surfaces in $M$) is determined by its vertical
bundle, which is an integrable rank $2$ distribution $\mathcal{H}$ on $F$, called the \textit{Lax distribution}.
Since the f\/ibre of $\mathcal{H}$ at $[\ell']\in PS'$ projects onto the $\alpha$-plane $S_{\pi_M
\ell'}\otimes\vspan{\ell'}$ in $T_{\pi_M\ell'} M$, it follows that the Lax distribution is given equivalently by an
inclusion $L_A\colon\pi_M^*S\to TF\otimes\mathcal{O}_F(1)$ called the \textit{Lax pair}, where $\mathcal{O}_F(-1)$ is
the f\/ibrewise tautological bundle of $PS'$.
I index the Lax pair $L_A$ to avoid confusion with the line bundle $L$ introduced previously: with respect to
a~trivialization of $\pi_M^*S$, $L_A$ really is a~pair $(L_0,L_1)$ of weighted vector f\/ields on $F$, and the
integrability of $\mathcal{H}$ means equivalently that $[L_0,L_1]\in\vspan{\{L_0,L_1\}}$.

The Lax pair can be made concrete by choosing trivialization $\lambda\colon F\to
{\mathbb{C}}{\rm P}^1={\mathbb{C}}\cup\{\infty\}$ of $F\to M$ (i.e., $(\pi_M,\lambda)$ identif\/ies $F$ with
$M\times {\mathbb{C}}{\rm P}^1$), which is essentially the same thing as a~trivialization of $S'$.
Together with a~trivialization of $S$, this gives
\begin{gather*}
L_0 = X_{00'} + \lambda X_{01'} + f_0(\lambda) \partial_\lambda,
\qquad
L_1 = X_{10'} + \lambda X_{11'} + f_1(\lambda) \partial_\lambda,
\end{gather*}
where $f_0$ and $f_1$ are functions on $F$ which are cubic in $\lambda$, $\partial_\lambda$ is the standard vector f\/ield
tangent to the ${\mathbb{C}}{\rm P}^1$-f\/ibres, and $X_{AA'}$ are pullbacks to $F$ of a~frame on $M$ corresponding
to the trivialization of $TM\cong S\otimes S'$.
This provides a~representative metric
\begin{gather*}
g= \theta_{00'}\theta_{11'}-\theta_{01'}\theta_{10'}
\end{gather*}
for the conformal structure, where $\theta_{AA'}$ is the dual coframe.
There is of course now the freedom to change the trivializations of $S$ and $S'$: then $\lambda$ undergoes a~M\"obius
transformation whose coef\/f\/icients are functions on $M$, while $L_0$ and $L_1$ get replaced by independent linear
combinations, again with functions on $M$ as coef\/f\/icients.

Now suppose $\ell\colon L\to S$ is a~selfdual $\beta$-surface foliation.
Then, pulling back to $F$ and composing with Lax pair yields a~line subbundle $L_A(\ell)$ of the Lax distribution which
is also tangent to the inverse image of the curves in $Z$ (which project to the $\beta$-surface foliation in $M$).
Let $(x_0,x_1)$ be components of the (local) leaf space projection $M\to N$ of the $\beta$-surface foliation, and let
$(\partial_{x_0}, \partial_{x_1})$ be lifts of the induced coordinate vector f\/ields on $N$.

Then the trivializations of $S$ and $S'$ can be chosen so that $L_0$ is a~section of $L_A(\ell)=L_A\ell^A$, that it is
tangent to the f\/ibres of $\lambda$, and that $L_1$ projects to the coordinate vector f\/ields on $N$ at $\lambda=0$ and
$\lambda=\infty$.
This means that the Lax pair takes the form
\begin{gather}\label{lax}
L_0 = \phi_0 + \lambda \phi_1,
\qquad
L_1 = \partial_{x_0}+\alpha_0 + \lambda (\partial_{x_1}+\alpha_1) + a(\lambda) \partial_\lambda,
\end{gather}
where $\phi_0$, $\phi_1$, $\alpha_0$ and $\alpha_1$ are vector f\/ields tangent to the f\/ibres of $M\to N$ (the
$\beta$-surface foliation), pulled back to $F$.
Since the Lax pair preserves the lifted $\beta$-surface foliation, it preserves the $3$-form $dx_0\wedge dx_1\wedge
d\lambda$ up to scale, and hence
\begin{gather*}
a(\lambda)=a_0(x_0,x_1)+a_1(x_0,x_1)\lambda+a_2(x_0,x_1)\lambda^2+a_3(x_0,x_1)\lambda^3
\end{gather*}
is a~function of $\lambda$, $x_0$, $x_1$ only.
The integrability of the Lax pair now reduces to
\begin{gather*}
[L_0,L_1]=b(\lambda)L_0=\big(b_0+b_1\lambda+b_2\lambda^2\big)L_0
\end{gather*}
for some functions $b_0$, $b_1$, $b_2$ on $M$.

$L_1$ pushes forward to the correspondence space between $N$ and $Y$ and this gives the projective spray on
${\rm P}(TN)\cong N\times{\mathbb{C}}{\rm P}^1$:
\begin{gather*}
\partial_{x_0} + \lambda \partial_{x_1} +
\bigl(a_0(x_0,x_1)+a_1(x_0,x_1)\lambda+a_2(x_0,x_1)\lambda^2+a_3(x_0,x_1)\lambda^3\bigr) \partial_\lambda,
\end{gather*}
where $x_0$, $x_1$, $\partial_{x_0}$, $\partial_{x_1}$ are now simply coordinates and coordinate vector f\/ields on~$N$.
The coef\/f\/icients~$a_j$ are related to the connection coef\/f\/icients~$\Gamma^i_{jk}$ of a~compatible projective connection
(where $\Gamma^i_{jk}=\Gamma^i_{kj}$) in the given trivialization of ${\rm P}(TN)$.
In fact~\cite{Hit:gse},
\begin{gather*}
a_0=\Gamma^1_{00},
\qquad
a_1=2\Gamma^1_{01}-\Gamma^0_{00},
\qquad
a_2=-2\Gamma^0_{01}+\Gamma^1_{11},
\qquad
a_3=-\Gamma^0_{11}.
\end{gather*}
If one now sets $\gamma_0=\Gamma^0_{00}+\Gamma^1_{01}$, $\gamma_1=\Gamma^0_{10}+\Gamma^1_{11}$,
$b(\lambda)=c(\lambda)-\frac13 a'(\lambda)$, where $c(\lambda)=c_0+c_1\lambda+c_2\lambda^2$ for functions $c_0$, $c_1$,
$c_2$, then the equation
\begin{gather*}
0=[L_1,L_0]+b(\lambda)L_0= [\partial_{x_0}+\alpha_0 + \lambda (\partial_{x_1}+\alpha_1),\phi_0 + \lambda \phi_1]
+a(\lambda)\phi_1+b(\lambda) (\phi_0 + \lambda \phi_1)
\end{gather*}
reduces to the condition $c_2=0$ and the following system:
\begin{gather*}
\partial_{x_0}\phi_0 + [\alpha_0,\phi_0] + \big(c_0-\tfrac23\gamma_0\big)\phi_0 + \Gamma^0_{00}\phi_0 + \Gamma^1_{00}\phi_1 = 0,
\\
\partial_{x_0}\phi_1 + [\alpha_0,\phi_1] + \big(c_0-\tfrac23\gamma_0\big)\phi_1 + \Gamma^0_{01}\phi_0 + \Gamma^1_{01}\phi_1
\\
\qquad
= -(\partial_{x_1}\phi_0 + [\alpha_1,\phi_0] + \big(c_1-\tfrac23\gamma_1\big)\phi_0 + \Gamma^0_{10}\phi_0 + \Gamma^1_{10}
\phi_1),
\\
\partial_{x_1}\phi_1 + [\alpha_1,\phi_1] + \big(c_1-\tfrac23\gamma_1\big)\phi_1 + \Gamma^0_{11}\phi_0 + \Gamma^1_{11}\phi_1 = 0.
\end{gather*}
The freedom left in the trivialization of $S$ takes the form $L_0\mapsto pL_0$ and $L_1\mapsto L_1+q L_0$ for functions
$p$, $q$ with $p$ nonvanishing, which changes $c(\lambda)$ by $L_0(q)-L_1(p)=\phi_0q-(\partial_{x_0}+\alpha_0)p
+\lambda(\phi_1q-(\partial_{x_1}+\alpha_1)p)$.
This can be used to set $c_0=c_1=0$, in which case the remaining interesting freedom (together with the coordinate
freedom in $x_0$, $x_1$ and the gauge freedom in the lifts of~$\partial_{x_0}$,~$\partial_{x_1}$) is just $L_1\mapsto L_1+q
L_0$ with $q$ a~function of $(x_0,x_1)$.

In this situation, there is a~natural interpretation of the above system on $N$ as a~gauge-theoretic system with gauge
group $\Diff_{2}$, the group of dif\/feomorphisms of a~surface (with Lie algebra $\mathop{\mathfrak{diff}}_2$, the
vector f\/ields on the surface): $\nabla^\alpha={\rm d}+\alpha=(\partial_{x_0}+\alpha_0,\partial_{x_1}+\alpha_1)$
def\/ines a~connection on a~bundle of surfaces over $N$, while $\phi=(\phi_0,\phi_1)$ is a~weighted $1$-form on $N$ with
values in the Lie algebra bundle ${\underline{\smash{\mathop{\mathfrak{diff}}}}}_2$ of vertical tangent vector f\/ields.
If $\mathcal{O}_N(1)$ denotes a~cube root of $\wedge^2TN$, then $\phi$ can be viewed as a~section of
$\mathcal{O}_N(2)\otimes T^{*} N\otimes{\underline{\smash{\mathop{\mathfrak{diff}}}}}_2$.
The gauge-theoretic system then makes $(\nabla^\alpha,\phi)$ into a~special case (albeit inf\/inite dimensional) of the
following projectively-invariant object.

\begin{defn}
Let $N$ be a~projective surface, over which there is a~principal $G$-bundle $P$ with adjoint bundle
$\underline{\mathfrak{g}}:=P\times_G\mathfrak{g}$.
Then a~pair $(\nabla^\alpha,\phi)$, where $\nabla^\alpha$ is a~$G$-connection and $\phi$ is a~section of
$\mathcal{O}_N(2)\otimes T^{*} N\otimes\underline{\mathfrak{g}}$ will be called a~\emph{projective pair} if\/f
\begin{gather*}
D^{\nabla^\alpha}\phi = \tfrac12{\rm d}^{D,\nabla^\alpha}\phi,
\end{gather*}
where $D$ is any connection in the projective class.
\end{defn}
In the projectively f\/lat case, this equation is a~reduction of the selfdual Yang--Mills equation on ${\mathbb{R}}^{2,2}$
or ${\mathbb{C}}^4$ by two null translations along a~$\beta$-plane, denoted $H_{\mathit{ASD}}$ in~\cite{MaWo:ist}
(cf.~\cite{Taf:2dr,TaWo:nkv}), although the projective invariance is not noted there.
(The fact that reductions of the selfduality equations for conformal structures result in gauge f\/ield equations given by
reductions of selfdual Yang--Mills equations is a~general phenomenon, which I have tried to systematize in my work on
integrable background geometries~\cite{Cal:ibg}, and that formalism gives another way to obtain some of the
results of the present paper.)

\begin{thm}
Let $(\nabla^\alpha,\phi)$ be a~projective pair, on a~projective surface $(N,[D])$, whose gauge group is a~subgroup of
$\Diff_{2}$ and such that the image of $\phi$ spans each tangent space of the associated surface bundle $M$.
Then $M$ is naturally equipped with a~selfdual conformal $4$-structure for which $M\to N$ is a~selfdual $\beta$-surface
foliation.
\end{thm}

\begin{proof}
This is simply a~matter of reversing the above arguments.
The pullback of $M\to N$ to ${\rm P}(TN)$ is a~$5$-manifold $F$ f\/ibering over $M$ with ${\mathbb{C}}{\rm P}^1$
f\/ibres.
The pullback of $\nabla^\alpha$ to $F$ gives a~lift of the projective spray to line subbundle $\vspan{L_1}\subset TF$,
$\phi$ def\/ines another line subbundle~$\vspan{L_0}$ and the projective pair equation on $(\nabla^\alpha,\phi)$ is
equivalent to the involutivity of $\vspan{\{L_0,L_1\}}$.
By construction~$L_0$ and $L_1$ have the form~\eqref{lax} in a~suitable trivialization.
Hence the associated frame on $M$ determines a~selfdual conformal structure~\cite{Pen:nlg}, and it is easy to see that
the f\/ibres of $M\to N$ form a~$\beta$-surface foliation.
By construction, there is a~lift of this foliation to the correspondence space $F$ which is Lie-preserved by the Lax
pair.
Hence there is a~Lie lift of the Lax pair to the area line bundle over the (lifted) foliation.
This def\/ines the selfdual connection realizing $M\to N$ as a~selfdual $\beta$-surface foliation.
(This is essentially the Penrose--Ward transform in the twistor theory for selfdual $\beta$-surface foliations which I
have already enunciated.)
\end{proof}

In the complex and real analytic categories, every selfdual space with a~selfdual $\beta$-surface foliation arises
locally in this way.
To check the same holds in the smooth category, it is enough to show that any selfdual $\beta$-surface foliation has
a~lift to $F$ of the form above.

Let $\ell\colon L \to S$ def\/ine a~$\beta$-surface foliation of $M$ with canonical connection $\nabla$.
Then $\mathcal{T}^\nabla\ell=0$, which implies that for any compatible metric with Levi-Civita connection $D$,
$D^\nabla\ell=\delta\otimes\varepsilon$, where $\varepsilon$ is the skew bilinear form on $S$ induced by the metric and
$\delta$ is an antiselfdual spinor f\/ield.
(In indices, $D^\nabla_{A'A}\ell_B=\delta_{A'}\varepsilon_{AB}$.) A tangent vector f\/ield $X=\beta\otimes\ell$ to $\ell$
is then lifted to $TF$ by taking, for each $[\pi]\in F$, the $D$-horizontal lift at $[\pi]$ (sometimes written
$\beta^{A'}\ell^A D_{A'A}$) minus $\ip{\beta\otimes\delta,\pi\otimes\pi}$ (i.e.,
$\beta^{A'}\delta^{B'}\pi_{A'}\pi_{B'}$), which is a~vertical vector f\/ield because it has homogeneity two in $\pi$.
I leave it to the enthusiastic reader to check that these lifts are Lie-preserved along the Lax pair (using the fact
that $\nabla$ is trivial along $\alpha$-surfaces).
This is all that one needs for all the previous arguments.

There is a~gauge-theoretic interpretation of the selfdual connection on $L\hookrightarrow S$: using this inclusion,
$L\otimes L$ is the subbundle of $\wedge^2_+TM$ corresponding to the area line bundle of the foliation, and the
connection is induced by the Lie lift of the Lax pair to the pullback of $L\otimes L$ to $F$.
In particular, the connection is f\/lat if\/f (locally) there is a~Lie-preserved area form along the lifted $\beta$-surface
foliation, i.e., if\/f the gauge group of $(\nabla^\alpha,\phi)$ reduces to $\SDiff_{2}$, the group of area
preserving dif\/feomorphisms.

\section{The Dunajski--West construction}

The next goal of this paper is to relate the preceding theory to the results of Dunajski and West for selfdual conformal
$4$-manifolds with null conformal vector f\/ields.
There are two subtasks here: f\/irst to show that their situation is a~special case of the present one; second to show how
their classif\/ication f\/its into this framework.
As mentioned at the start of this paper, the key to the Dunajski--West classif\/ication is the following observation.

\begin{lem}[\cite{DuWe:acs}]
\label{int}
Let $M$ be a~conformal $4$-manifold with a~$($nonzero$)$ null conformal Killing vector~$K$.
Then the integral curves of $K$ are geodesics, and the selfdual and antiselfdual null $2$-plane distributions in
$K^\perp$ are both integrable.
\end{lem}

\begin{proof}
I give a~tensorial proof of this result, referring to~\cite{DuWe:acs} for the spinorial approach, which is, of course,
very relevant for the present paper.
The geodesic property is obvious (and is a~general result about Killing vector f\/ields of constant length): the conformal
Killing equation means that for any compatible metric $g$, $D^gK$ is a~section of $\mathop{\mathfrak{co}}(TM)$
(i.e., skew plus a~multiple of the identity); hence for any vector f\/ield $X$ orthogonal to $K$, $0=\ip{D^g_X K,
K}=-\ip{D^g_K K,X}$, so $D^g_KK$ is a~multiple of $K$.
The second assertion amounts to the fact that if $X$ is a~\emph{null} vector f\/ield orthogonal to $K$, then
$\vspan{\{K,X\}}$ is an integrable distribution on the open set where it has rank $2$.
Since this distribution is totally null, $[X,K]$ lies in the distribution if and only if it is orthogonal to it! Now
$g(D^g_K X-D^g_X K,K)=g(D^g_K X,K)= -g(X,D^g_K K)=0$ and $g(D^g_K X-D^g_X K,X)=g(D^g_X K,X)=0$ (by the conformal Killing
equation and the fact that $X$ is null).
\end{proof}

For the f\/irst subtask, let $M$ be selfdual with a~null conformal vector f\/ield.
Then $K$ def\/ines a~$\beta$-surface foliation $\ell$ and an $\alpha$-surface foliation $\ell'$ on the open set where it
is nonzero.
I claim that the $\beta$-surface foliation is selfdual and have two ways to see this, one twistorial, the other by
direct computation.
The twistorial approach amounts to observing that~$K$ def\/ines a~vector f\/ield~$X$ on the twistor space~$Z$ transverse to
twistor lines and then showing that the span of~$X$ can be continued as a~distribution through the divisor $\mathcal{D}$
on which it vanishes (which is the divisor corresponding to the $\alpha$-surface foliation to which $K$ is tangent):
Dunajski and West do this in~\cite[Lemma 3]{DuWe:acs}.
For a~more abstract argument, let $s$ be a~section of the degree~$1$ divisor line bundle
$[\mathcal{D}]=\mathcal{L}\otimes\mathcal{O}_Z(1)$ of $\mathcal{D}$ vanishing nondegenerately on $\mathcal{D}$.
Then $\hat\ell:=s^{-1}X$ is a~meromorphic section of $\Hom(\mathcal{L}\otimes\mathcal{O}_Z(1), TZ)$.
However, for any twistor line~$\hat x$, the projection of~$\hat\ell$ onto the normal bundle of~$\hat x$ is a~meromorphic
section of a~trivial bundle which is bounded on~$\hat x\setminus \mathcal{D}$, hence extends to a~trivialization over~$\hat x$ by Liouville's theorem.
Since this is true for any twistor line, $\hat\ell$ is actually holomorphic.
(A similar continuation result was used in~\cite{CaPe:sdc} in the context of geodesic symmetries in Einstein--Weyl
geometry.)

The direct computation is also easy: $\alpha$-surface foliations also have a~canonical connection by the same argument
as for $\beta$-surface foliations, but on a~selfdual conformal manifold, such connections are always selfdual by an easy
computation (in fact they correspond to the `degree zero part'~$\mathcal{L}$ above of the divisor line bundle
$[\mathcal{D}]$); however, since the conformal vector f\/ield~$K$ is essentially~$\ell\otimes\ell'$, the two canonical
connections are (up to gauge equivalence) dual to each other.

\begin{rem}
As an aside for the occasional reader with an interest in selfdual conformal submersions~\cite{Cal:sde}, I note that the
theory given there generalizes as follows.
A line subbundle $\xi\colon L\to TM$ of the tangent bundle will called \textit{semi-conformal submersion} if there is
a~connection $\nabla$ on $L$ such that $\sym_0D^\nabla\xi=0$ (vanishing symmetric tracefree part, or, in indices,
$(D^\nabla\xi)_{(A'B')(AB)}=0$): clearly $\nabla$ is uniquely determined by this condition and will be called the
canonical connection of $\xi$.
A semi-conformal submersion will be called \textit{selfdual} if\/f $\nabla$ has selfdual curvature.

If $\xi$ is non-null, then after identifying $L$ with the weight~$1$ density line bundle of~$M$ so that~$\xi$ has unit
length, it becomes a~conformal submersion in the sense of~\cite{Cal:sde}.
On the other hand, in the null case, Lemma~\ref{int} generalizes straightforwardly to this context, with the result that
$\xi=\ell\otimes\ell'$ for an $\alpha$-surface foliation~$\ell'$ and a~$\beta$-surface foliation~$\ell$.
The canonical connection of~$\xi$ is easily seen to be the tensor product of the canonical connections of $\ell$ and~$\ell'$.
Since the latter is always selfdual, $\xi$ is selfdual if\/f~$\ell$ is.
\end{rem}

Now to the second subtask: the Dunajski--West classif\/ication.
The observation here is a~simple one: a~selfdual conformal manifold, with a~$\beta$-surface foliation constructed from
a~projective pair, has a~(null) conformal vector f\/ield tangent to foliation if and only if the gauge group of the
projective pair reduces to the group of dif\/feomorphisms of a~surface which preserve a~given vector f\/ield.
This follows immediately from the bijective nature of the correspondence I have given (it is therefore functorial up to
gauge equivalence).
Choose coordinates $(t,z)$ on the surfaces so that the preserved vector f\/ield is $\partial_t$.
This realizes the gauge group (locally) as a~central extension
$\mathcal{O}^{\times}\rtimes\Diff_{1}$ of the group of dif\/feomorphisms of a~$1$-manifold by the
space of nonvanishing functions on that manifold: the Lie algebra $\mathfrak o\rtimes\mathop{\mathfrak{diff}}_1$
consists of vector f\/ields of the form $f(z)\partial_t+g(z)\partial_z$.
For easier comparison with~\cite{DuWe:acs}, I will write $(x,y)$ for $(x_0,x_1)$.
The Lax pair~\eqref{lax} now takes the form
\begin{gather*}
L_0 = (p(z) + \lambda q(z))\partial_t+(u(z) + \lambda v(z))\partial_z,
\\
L_1 = \partial_{x} + \lambda \partial_{y} + (C(z)+\lambda D(z))\partial_t+ (E(z)+\lambda F(z))\partial_z + a(\lambda)
\partial_\lambda,
\end{gather*}
where all functions also depend on $(x,y)$ but I do not denote this explicitly.
The conformal structure is therefore represented by the metric
\begin{gather*}
g=\bigl({\rm d} t-C{\rm d} x-D{\rm d} y)(u{\rm d} x+v{\rm d} y)
+({\rm d} z-E{\rm d} x-F{\rm d} y)(p{\rm d} x+q{\rm d} y)
\end{gather*}
and the \textit{twist} of $K=\partial_t$ (i.e., $*(\alpha\wedge {\rm d}\alpha)$ where
$\alpha=g(K,\cdot)$) is $v u_z-u v_z$ up to sign.
Notice that the twist vanishes if and only if $u/v$ is independent of $z$.
Now $\lambda=u/v$ is the hypersurface on which $\partial_t$ is tangent to the Lax distribution, which projects to the
divisor of the $\alpha$-surface foliation in the twistor space $Z$.
Thus $K$ is twist-free if and only if the divisor of the $\alpha$-surface foliation is the pullback of a~divisor in $Y$
(which corresponds to a~geodesic congruence in $N$).

There is a~natural gauge-theoretic approach to solve the Lax equation with this gauge group:
f\/irst solve the projective pair equation for $\Diff_{1}$;
then the full solution on the central extension is given by quadratures.
The $\Diff_{1}$ part of the Lax pair is
\begin{gather*}
L_0' = (u(z) + \lambda v(z))\partial_z,
\qquad
L_1' = \partial_{x} + \lambda \partial_{y} + (E(z)+\lambda F(z))\partial_z + a(\lambda) \partial_\lambda,
\end{gather*}
and I shall suppose $uv\neq 0$, so that, without loss, $v\neq 0$.
Since the $\Diff_{1}$ projective pair equation is nonlinear, it is not obvious that one can solve it.
However, the divisor $\lambda=u/v$ appearing in the twist-free case above is already present at the $\Diff_{1}$ level,
and it turns out that even when~$u/v$ depends on~$z$, there are still divisors in~$Y$ corresponding to geodesic
congruences in~$N$, except that they too depend on~$z$.

To see this explicitly, I f\/irst make a~a coordinate change $\tilde z=f(x,y,z)$ (with $f_z\neq 0$) and a~change of
trivialization $\smash{\tilde L}_1'=L_1'+g(x,y,z) L_0'$, $\smash{\tilde L}_0'=h(x,y,z)L_0'$ to set $E=F=0$ and $v=1$, so
that writing $u=-\beta$ and dropping tildes yields a~simpler form
\begin{gather*}
L_0' = (\lambda - \beta(z))\partial_z,
\qquad
L_1' = \partial_{x} + \lambda \partial_{y} + a(\lambda) \partial_\lambda
\end{gather*}
for the Lax pair.
However, the change of trivialization comes at a~price: $[L_0',L_1']=b(z,\lambda)L_0'$, where $b$ now depends on $z$ in
general.
Despite this, the Lax equation is very simple:
\begin{gather*}
b=\beta_y-\frac{a(\lambda)-a(\beta)}{\lambda-\beta},
\qquad
\beta_x+\beta\beta_y -a(\beta)=0.
\end{gather*}
The f\/irst equation determines $b$ (as a~quadratic in $\lambda$), while the second means that for each f\/ixed $z$,
$\lambda=\beta(x,y,z)$ is a~section of $P(TN)$ tangent to the geodesic spray.
Hence, as claimed above, $z$ parameterizes a~family of geodesic congruences in $N$ (corresponding to divisors in $Y$
transverse to minitwistor lines) generalizing the twist-free case (where $\beta_z=0$).

I end this section with a~gauge-theoretic characterization of the twist-free condition.
For this, I~abandon the cavalier changes of trivialization made above, and insist instead that $b(\lambda)=-\frac13
a'(\lambda)$ (independent of $z$).
However, the coordinate change $\tilde z=f(x,y,z)$ can still be applied to set $v=1$, and with $u=-\beta$, the
twist-free case is still $\beta_z=0$.
Assuming this, the extra ingredients~$E$ and~$F$ satisfy $E_z=\beta_y+\frac23 a'(\beta)$ and $F_z=\frac16 a''(\beta)$,
so that $E$ and $F$ are af\/f\/ine in~$z$.
By translation of $z$, and the addition of a~($z$-independent) multiple of~$L_0'$ to~$L_1'$, any projective pair which
is af\/f\/ine in~$z$, with~$L_0'$ nonvanishing but independent of~$z$, can be supposed to have~$E$ and~$F$ linear in~$z$.
Hence the Lax pair reduces to
\begin{gather}
\label{geom-lax}
L_0' = (\lambda - \beta)\partial_z,
\qquad
L_1' = \partial_{x} + \lambda \partial_{y} +\bigl(\beta_y+\tfrac23 a'(\beta)+\tfrac16 (\lambda-\beta)
a''(\beta)\bigr)z\partial_z +a(\lambda) \partial_\lambda
\end{gather}
and $[L_0',L_1']=-\frac13 a'(\lambda) L_0$ if and only if $\beta_x+\beta\beta_y -a(\beta)=0$, i.e.,
$\lambda=\beta$ def\/ines a~geodesic congruence on $N$.

\begin{thm}
Let $M$ be a~selfdual conformal manifold with a~null conformal vector field $K$ constructed from an
$\mathcal{O}^{\times}\rtimes \Diff_{1}$ projective pair $(\nabla^\alpha,\phi)$ on a~projective surface $N$.
Then $K$ is twist-free if and only if the induced $\Diff_{1}$ projective pair reduces to the affine group $\Aff_{1}$
of a~line, with $\phi$ purely translational.

Furthermore, such $\Aff_{1}$ projective pairs correspond bijectively $($locally, up to gauge and coordinate
transformations, and $\nabla^\alpha\mapsto\nabla^\alpha+f(x,y)\phi)$ to geodesic congruences on~$N$.
\end{thm}

Following the methods of the appendix, it is easy to see that the linear ($\GL_1$) part of the $\Aff_{1}$ connection
$\nabla^\alpha$ is the Penrose--Ward transform of $\mathcal{L}_\mathcal{C}$, where $\mathcal{C}$ is the divisor
corresponding to the geodesic congruence, with divisor line bundle
$[\mathcal{C}]=\mathcal{L}_{\mathcal{C}}\otimes\mathcal{O}_Y(1)$.
On the other hand, $\phi$ is the (coupled) Penrose transform of the section of $[\mathcal{C}]$ (unique up to scale)
vanishing on $\mathcal{C}$.
In particular $\mathcal{L}_\mathcal{C}$ is trivial (and so $-3\mathcal{C}$ is a~canonical divisor) if $\beta_y+\tfrac23
a'(\beta)=0$ and $a''(\beta)=0$.
This will be useful in the f\/inal section.

\begin{rem}
Dunajski and West proceed in a~slightly dif\/ferent way.
They use the freedom to change the trivializations and the $(t,z)$ coordinates to set $p=1$, $q=0$ and $v=1$, leaving,
in $L_0$, only one nontrivial function $\beta:=-u$.
They then consider separately the twisting case (i.e., $\beta_z\neq 0$) and twist-free case (i.e.,
$\beta_z=0$), and reduce the projective pair equation to quadratures.
The two cases can be considered together by f\/ixing the gauge dif\/ferently in the twisting case.
In general, one can change variables so that $\beta(x,y,z)=\gamma(x,y)+cz$, with $c=0$ in the twist-free case.
Rescaling $L_0$ then gives $L_0= H(z)\partial_t +(\lambda-\beta(z))\partial_z$, for some function $H$.
The involutivity condition $[L_0,L_1]\in\vspan{L_0}$ may then be reduced to: $E=(a_1+\gamma a_2+\gamma^2 a_3-\gamma_y)z$
and $F= a_3z (\gamma+cz)+(a_2+\gamma a_3)z$ (where~\cite{DuWe:acs} would have the simpler formulae $E=a(z)$ and $F=0$ in
the twisting case); $C_z=0$ and $D_z=-a_3 H$ (where~\cite{DuWe:acs} has a~more complicated formula for $C_z$ in the
twisting case); and
\begin{gather*}
H_x+\beta H_y + (E+\beta F)H_z =0,
\qquad
\gamma_x +\gamma \gamma_y - a(\gamma) = 0.
\end{gather*}
Setting, $C=0$ and $D=-a_3 G$, where $G_z=H$, these may be solved by quadratures, given a~geodesic congruence on $N$.
\end{rem}

\section{Hyperk\"ahler, hypercomplex,\\ and scalar-f\/lat K\"ahler structures}

I end this paper by asking one of the questions raised by~\cite{DuWe:acs}: when does a~selfdual conformal manifold $M$
constructed in this way admit a~compatible Einstein metric? In the Ricci-f\/lat case, this means $M$ is (locally)
hyperk\"ahler in the complexif\/ied or kleinian signature sense (sometimes called pseudo-hyperk\"ahler).
This leads to a~more general question concerning the existence (with my orientation convention) of antiselfdual complex
structures (including hypercomplex structures and scalar-f\/lat K\"ahler metrics) or antiselfdual integrable involutions
or $\alpha$-surface foliations (including pseudo-hypercomplex structures and selfdual pseudo-K\"ahler or null-K\"ahler
metrics).

There is a~standard approach to this question in the context of reductions, using the theory of divisors (see
also~\cite{Cal:sde,Cal:ibg,CaPe:sdc,Dun:4ps}).
An antiselfdual (null- or pseudo-) complex structure $J$ corresponds to a~degree two divisor $\mathcal{D}$ in the
twistor space $Z$, which locally has the form $\mathcal{D}=\mathcal{D}_1+\mathcal{D}_2$ with
$\mathcal{D}_1=\mathcal{D}_2$ in the null case, but disjoint otherwise.
It is then natural to suppose that a~selfdual $\beta$-surface foliation is compatible with $\mathcal{D}$ in the sense
that the corresponding foliation of $Z$ also foliates $\mathcal{D}$, so that the latter is $\pi^*\mathcal{C}$ for
a~degree two divisor $\mathcal{C}=\mathcal{C}_1+\mathcal{C}_2$ in the leaf space~$Y$, where $\pi\colon Z\to Y$ is the
quotient map.
I discuss the theory of degree two divisors in minitwistor spaces in an appendix: they correspond to a~pair of geodesic
congruences on a~projective surface~$N$, which coincide in the null case ($\mathcal{C}_1=\mathcal{C}_2$).

As is well known (see~\cite{CaPe:sdc,Dun:4ps} and references therein), $J$ is the complex structure (or involution) of
a~compatible scalar-f\/lat (null- or pseudo-) K\"ahler metric if\/f $\mathcal{D}$ is a~divisor for $\mathcal{O}_Z(2)$, the
square root of the anticanonical bundle of $Z$: in other words, $[\mathcal{D}]\otimes \mathcal{O}_Z(-2)$ is trivial.
Now, following the notation of Theorem~\ref{main}, the vertical bundle of $\pi$ is $\mathcal{L}\otimes\mathcal{O}_Z(1)$,
where $\mathcal{L}$ has degree zero, and so an exact sequence shows that $\mathcal{O}_Z(3)=\mathcal{L}\otimes
\pi^*\mathcal{O}_Y(3)$.
Thus $[\mathcal{D}]\otimes\mathcal{O}_Z(-2)=\mathcal{L}^{-2/3}\otimes\pi^*([\mathcal{C}]\otimes\mathcal{O}_Y(-2))$ which
is trivial if\/f $\mathcal{L}^{2/3}\cong\pi^*([\mathcal{C}]\otimes\mathcal{O}_Y(-2))$.
This implies a~reduction of gauge group for the projective pair $(\nabla^\alpha,\phi)$, since the induced connection on
the line bundle of f\/ibrewise area forms is basic relative to the f\/ibration $M\to N$: the gauge group thus reduces to the
group $\HDiff_{2}$ of dif\/feomorphisms preserving an area form up to scale, which is an extension
$\SDiff_{2}\rtimes\GL_1$ of $\GL_1$ (${\mathbb{R}}^{\times}$ or
${\mathbb{C}}^{\times}$) by the group $\SDiff_{2}$ of area preserving dif\/feomorphisms.
The associated $\GL_1$-connection is given by the Penrose--Ward transform of $\mathcal{L}\cong
\pi^*([\mathcal{C}]^{3/2}\otimes\mathcal{O}_Y(-3))$.

In the nonnull case, where $\mathcal{D}_1\neq \mathcal{D}_2$, $J$ is (locally) part of a~(complexif\/ied or kleinian)
hypercomplex structure if\/f $[\mathcal{D}_1-\mathcal{D}_2]$ is trivial, i.e., $\mathcal{D}_1-\mathcal{D}_2$ is
the divisor of a~meromorphic function (which is then locally a~f\/ibration of $Z$ over ${\mathbb{C}}{\rm P}^1$).
When $[\mathcal{D}_1-\mathcal{D}_2]= \pi^*[\mathcal{C}_1-\mathcal{C}_2]$, this actually forces
$[\mathcal{C}_1-\mathcal{C}_2]$ to be trivial already.

\begin{thm}\label{hk}
Let $N$ be a~projective surface whose minitwistor space $Y$ contains a~degree two divisor $\mathcal{C}$ transverse to
minitwistor lines.
Let $M\to N$ be a~$\beta$-surface foliation over $N$ induced by a~projective pair $(\nabla^\alpha,\phi)$ on $N$, and let
$\mathcal{D}$ be the pullback of $\mathcal{C}$ to the twistor space $Z\to Y$.
\begin{enumerate}\itemsep=0pt
\item[$(i)$]
$\mathcal{D}$ is the divisor of a~$($locally$)$ scalar-flat $($null- or pseudo-$)$ K\"ahler
structure if and only if the gauge group of $(\nabla^\alpha,\phi)$ reduces to $\HDiff_{2}$, $\phi$ takes values in
$\underline{\smash{\mathop{\mathfrak{sdiff}}}}_{2}$, and the induced $\GL_1=\HDiff_{2}/\SDiff_{2}$ connection is the
Penrose--Ward transform of $[\mathcal{C}]^{3/2}\otimes\mathcal{O}_Y(-3)$.
\item[$(ii)$]
If $\mathcal{C}$ is not twice a~degree one divisor, then the complex structure corresponding to $\mathcal{D}$ is part of
a~hypercomplex structure if and only if there is a~fibration $f\colon Y\to{\mathbb{C}}{\rm P}^1$ with
$\mathcal{C}\subseteq f^{-1}(\{0,\infty\})$.
\end{enumerate}
In particular, $\mathcal{C}$ induces a~$($locally$)$ hyperk\"ahler structure on $M$ if both these conditions
hold.
\end{thm}

\begin{proof}
There is essentially nothing left to prove in (i), so I turn to (ii).
Then $\mathcal{D}=\pi^*\mathcal{C}=\pi^*(\mathcal{C}_1+\mathcal{C}_2)$ induces one of the complex structures of
a~hypercomplex structure (where $\pi\colon Z\to Y$ is the projection) if\/f $\pi^*[\mathcal{C}_1-\mathcal{C}_2]$ is
trivial, i.e., it admits a~nonvanishing section $s$ (unique up to scale).
The derivative of any such $s$ along the f\/ibres of $\pi$ must vanish, since it def\/ines a~bundle map from the vertical
bundle of $Z\to Y$, which has degree one on $Z$, to $\pi^*[\mathcal{C}_1-\mathcal{C}_2]$, which has degree zero: thus
$s$ is constant on the f\/ibres of $\pi$.
Hence $\pi^*[\mathcal{C}_1-\mathcal{C}_2]$ is trivial if\/f $[\mathcal{C}_1-\mathcal{C}_2]$ is, in which case
$\mathcal{C}_1$ and $\mathcal{C}_2$ are the zero and pole divisors of a~meromorphic function $Y\to
{\mathbb{C}}{\rm P}^1$.
\end{proof}

Unfortunately, this theorem does not completely characterize the case that (twice) a~degree one divisor pulls back to
a~f\/ibre of a~meromorphic function, and so does not cover all constructions of hyperk\"ahler metrics with a~pullback
divisor in the hyperk\"ahler family.
For instance, if a~hyperk\"ahler metric has a~null (perhaps only homothetic) Killing vector $K$ which is not
triholomorphic, then the f\/ibration of $Z$ over ${\mathbb{C}}{\rm P}^1$ does not descend to $Y$, and at most one of
its f\/ibres does.
(If $K$ is twist-free, the divisor of its $\alpha$-surface foliation is the only pullback divisor in the hyperk\"ahler
family, whereas if $K$ has twist, it preserves two distinct divisors in the hyperk\"ahler family, but one of them is the
divisor of its $\alpha$-surface foliation, which, as already noted, is not a~pullback from $Y$ in the twisting case.)

In the context of the Dunajski--West construction, it is natural to study the projective pairs for which the gauge group
reduces to the subgroup $G$ of $\HDiff_{2}$ commuting with a~vector f\/ield which preserves the area form def\/ining
$\HDiff_{2}$.

\begin{rem}
If one is interested in selfdual Einstein metrics, then results of Pedersen--Tod~\cite{PeTo:emh} and
Przanowski~\cite{Prz:kvf} (see also~\cite{Tod:sut}) can be adapted to justify in part the reduction to $G$.
Indeed, if $g$ is a~selfdual Einstein metric with a~conformal vector f\/ield $K$, so that $\mathcal{L}_Kg = f g$ for some
function $f$, then~\cite[Proposition 3.2]{PeTo:emh} shows that one of the following three possibilities occurs: (i)~$\nabla^g f$ is the only principal null direction of the Weyl curvature $W$ (so $W$ has Petrov type~$N$); (ii)~$f$~is
constant (so~$K$ is a~homothetic vector f\/ield) and $g$ is hyperk\"ahler; or (iii)~$f=0$ so that~$K$ is a~Killing vector
f\/ield.
A homothetic vector f\/ield on a~hyperk\"ahler manifold preserves the Levi-Civita connection, hence at least one (possibly
null) complex structure in the hyperk\"ahler family.
On the other hand, by~\cite{Prz:kvf,Tod:sut}, a~selfdual Einstein metric with a~Killing vector f\/ield $K$ is conformal to
a~scalar-f\/lat (possibly null- or pseudo-) K\"ahler metric for which $K$ is a~holomorphic Killing vector f\/ield.
There remains the issue, mentioned already in the hyperk\"ahler case, that corresponding degree two divisor in $Z$ may
not be a~pullback from~$Y$: the reduction to~$G$ only covers the case that it is.
\end{rem}

In any case, it is easy to see that $G$ is $\mathcal{O}^\times\rtimes\Aff_{1}$, which is solvable.
Explicitly, if the symmetry is generated by a~vector f\/ield $\partial_t$ with momentum map $z$, the Lie algebra
$\mathfrak o\rtimes\mathop{\mathfrak{af\/f}}_1$ consists of vector f\/ields of the form $f(z)\partial_t + (uz+w)\partial_z$.
Now the projective pairs $(\nabla^\alpha,\phi)$ with this gauge group such that $\phi$ takes values in
${\underline{\smash{\mathop{\mathfrak{sdiff}}}}}_2$, as in Theorem~\ref{hk}, are exactly those arising in the twist-free
case ($\beta_z=0$) of the Dunajski--West construction.
which have already been associated to degree one divisors in $Y$ (geodesic congruences in $N$).

There is a~tricky point however: to obtain a~scalar-f\/lat (null- or pseudo-) K\"ahler metric, the induced $\GL_1$
connection must be the Penrose--Ward transform of $[\mathcal{C}]^{3/2}\otimes \mathcal{O}_Y(-3)$ for a~degree two
divisor $\mathcal{C}$, whereas it is naturally the Penrose--Ward transform of $[\mathcal{C}']\otimes\mathcal{O}(-1)$ for
a~degree one divisor $\mathcal{C}'$.
To reconcile these conditions, another degree one divisor $\mathcal{C}''$ is needed, such that
$[\mathcal{C}'']^2\otimes\mathcal{O}(-2)\cong [\mathcal{C}']\otimes\mathcal{O}(-1)$: then
$\mathcal{C}=\mathcal{C}'+\mathcal{C}''$.
The simplest way to achieve this is to suppose that $[\mathcal{C}']=\mathcal{O}(1)$ and take
$\mathcal{C}''=\mathcal{C}'$ so that $[\mathcal{C}]^{3/2}\otimes\mathcal{O}(-3)$ is trivial and the $\GL_1$ connection
is f\/lat.
As noted by Hitchin~\cite{Hit:cme}, this means the projective spray can be reduced to the form
$\partial_x+\lambda\partial_y+a\partial_\lambda$, with $a=a(x,y)$ independent of $\lambda$.
One way to see this is to choose $y$ to be constant along the geodesics of the congruence.
Then $\beta=0=a(0)$ in~\eqref{geom-lax}, which yields the generalized pp-wave metrics of~\cite[\S~6.2]{DuWe:acs}.
Here, however, the f\/latness of the $\GL_1$ connection can be used to set $a'(0)=0=a''(0)$.
Inverting $\lambda$ and exchanging~$(x,y)$ gives Hitchin's form of the spray.
The full Lax pair is then
\begin{gather*}
L_0=(p(z)+\lambda q(z))\partial_t + \partial_z,
\qquad
L_1=\partial_x+\lambda\partial_y+(C(z)+\lambda D(z))\partial_t+a\partial_\lambda,
\end{gather*}
and a~translation of $t$ can be used to eliminate $p$.
Now $[L_0,L_1]=0$ gives $C_z=a q$, $D_z=q_x$ and $q_y=0$.
Modulo changes of $z$ coordinate and trivialization, the general solution can be put into the form $q=1$, $C=az$,
$D_z=0$.
Thus
\begin{gather*}
L_0=\partial_z+\lambda\partial_t,
\qquad
L_1=\partial_x+az \partial_t +\lambda(\partial_y+c\partial_t)+a\partial_\lambda,
\end{gather*}
where $a$ and $c$ are functions of $(x,y)$ only.
The corresponding selfdual conformal structure then contains a~family of null-K\"ahler metrics
\begin{gather*}
g=f\bigl({\rm d} z{\rm d} y -({\rm d} t-az {\rm d} x- c{\rm d} y){\rm d}x\bigr),
\end{gather*}
where $f$ is an arbitrary function of $x$ and $z$: the null complex structure is $J={\rm d} z\otimes
\partial_t+{\rm d} x\otimes \partial_y$ with null K\"ahler form $\omega=f {\rm d} x\wedge {\rm d}
z$.
When $a=0$ and $f=1$, this is a~pp-wave hyperk\"ahler metric.
As observed in~\cite{DuWe:acs}, $g$ is also hyperk\"ahler for $c=0$ and $f=1/z^2$: the $\partial_\lambda$ term in the
Lax pair can be eliminated by a~change $\lambda\mapsto z\lambda-t$ of the spectral parameter $\lambda$.

One might hope to obtain other hyperk\"ahler metrics from $\mathcal{O}^\times\rtimes\Aff_{1}$ projective pairs over
projective surfaces with a~minitwistor space f\/ibering over ${\mathbb{C}}{\rm P}^1$, but the obvious construction
with $\mathcal{C}=f^{-1}(\{0,\infty\})$, where $f\colon Y\to {\mathbb{C}}{\rm P}^1$ is the f\/ibration, yields only
pp-waves.
Indeed, the Penrose--Ward transform of $[\mathcal{C}]\otimes\mathcal{O}_Y(-2)$ is a~connection $\nabla$ on
$\mathcal{O}_N(1)$ inducing a~projective connection $D$ whose Ricci curvature ${\rm r}^D$ is skew (see the
appendix), since it is a~constant multiple of the curvature $F^\nabla$ of $\nabla$.
The $\Aff_{1}$ projective pair is then obtained by solving $\sym D^\nabla\phi=0$, but this has an integrability
condition: $0={\rm d}^DF^\nabla\sim {\rm d}^D{\rm r}^D$ in $\wedge^2T^{*} N\otimes T^{*} N$.
Hence the projective structure is f\/lat, $Y$ is an open subset of ${\mathbb{C}}{\rm P}^2$ and
$[\mathcal{C}]\otimes\mathcal{O}_Y(-2)$ is trivial ($f$ is then the restriction to $Y$ of a~standard f\/ibration
${\mathbb{C}}{\rm P}^2\setminus\{x\}\cong\mathcal{O}(1)\to {\mathbb{C}}{\rm P}^1$).
Thus $\nabla$ is f\/lat, which is the construction already considered.

\appendix \section{Minitwistor theory of projective surfaces}

In this appendix, I sketch some aspects of the minitwistor theory of projective surfaces, which are straightforward, but
not detailed in the existing literature.
Let $N$ be a~geodesically convex complex projective surface and $Y$ its space of geodesics.
The correspondence space is $P(TN)$ endowed with the projective spray.
Suppose that the anticanonical bundle $\wedge^2TN$ has a~cube root $\mathcal{O}_N(1)$ and let $\varepsilon$ denote the
natural section of $\wedge^2T^*N \mathcal{O}_N(3)$.
Here, and in the following, I~omit the tensor product sign when tensoring with a~line bundle.

It will be useful to introduce an algebraic bracket $T^{*} N\otimes TN \to \mathop{\mathfrak{gl}}(TN)$ def\/ined by
$\abrack{\gamma,X}=\gamma(X){\rm id} + \gamma\otimes X$: then projectively equivalent connections~$D$, $\tilde D$
dif\/fer by the $\mathop{\mathfrak{gl}}(TN)$-valued $1$-form $\abrack{\gamma,\cdot}$, which acts naturally on any tensor
representation.
The curvature $R^D$ of a~connection~$D$ in the projective class is necessarily of Ricci type, i.e.,
$R^D=\abrack{{\rm r}^D\wedge{\rm id}}$ for a~bilinear form ${\rm r}^D$ on~$TN$, and
${\rm r}^{D+\abrack{\gamma,\cdot}}={\rm r}^D+D\gamma-\gamma\otimes\gamma$.
The curvature of the projective structure is the projectively invariant section $C={\rm d}^D{\rm r}^D$ of
$\wedge^2T^{*} N\otimes T^{*} N$, which is an analogue of the Cotton--York tensor on a~conformal $3$-manifold.

The f\/irst observation is that the equation $D_XX=0$ for geodesic vector f\/ields is projectively invariant if $X$ is
a~section of $\mathcal{O}_N(-2) TN\cong \mathcal{O}_N(1)  T^{*} N$, because $\abrack{\gamma,X}\cdot \alpha=
\alpha(X)\gamma$ on $\mathcal{O}_N(1)  T^{*} N$: hence if $\alpha=\varepsilon(X,\cdot)$ and $\tilde
D=D+\abrack{\gamma,\cdot}$, then $\tilde D_X\alpha=D_X\alpha$.
This def\/ines a~lift of the geodesic spray to $\mathcal{O}_N(-2) TN$, which is a~line bundle over the correspondence
space $P(TN)$.
The f\/low is linear on f\/ibres and hence $\mathcal{O}_N(-2) TN$ is the pullback of a~degree minus one line bundle on $Y$,
and one easily sees that this is $\mathcal{O}_Y(-1)$, a~cube root of the canonical bundle $K_Y$.

The connection coef\/f\/icients $\hat\Gamma^A_{BC}$ on $\mathcal{O}_N(-2) TN$ are related to the Christof\/fel symbols
$\Gamma^A_{BC}$ in a~frame $X_A$ by $\hat\Gamma^A_{BC} = \Gamma^A_{BC}-\frac23\delta^A_C\Gamma^E_{BE}$.
The lifted spray may be written in indices as $\pi^AX_A + \pi^B\pi^C\hat\Gamma^A_{BC}\partial_{\pi^A}$ and its
projective invariance (under $\Gamma^A_{BC}\mapsto \Gamma^A_{BC}+\gamma_B\delta^A_C+\gamma_C\delta^A_B$) is easy to
check directly.

\subsection*{Penrose transforms}

The Penrose transform relates the cohomology of vector bundles over $Y$ to spaces of sections over $N$ in (usually)
kernels or cokernels (or quotients thereof) of natural dif\/ferential operators.
For instance, using the lifted spray above, the space $H^0(Y,\mathcal{O}_Y(1))$ of holomorphic sections of~$\mathcal{O}_Y(1)$ is readily identif\/ied with the space of sections $\phi$ of $\mathcal{O}_N(-1) TN$ satisfying the
abelian projective pair equation $D\phi=\frac12{\rm d}^D\phi$.

I shall now wheel out some machinery to compute the Penrose transform for a~larger class of bundles.
In doing so, I cheerfully accept the accusation of using a~sledgehammer, indeed two sledgehammers, to crack a~nut, but
they are, in my opinion, very appealing sledgehammers and have not been applied in this context before.

Projective surfaces are examples of parabolic geometries, which are Cartan geometries modelled on generalized f\/lag
varieties $G/P$, where $G$ is a~semisimple Lie group and $P$ is a~parabolic subgroup (i.e., the Lie algebra of
$P$ is the Killing perp of its maximal nilpotent ideal).
A Cartan geometry modelled on such a~homogeneous space is a~manifold~$M$ of the same dimension as~$G/P$, together with
a~principal $P$-bundle $\mathcal{G}\to M$ and a~principal $G$-connection $\omega$ on $\mathcal{G}\times_P G$ whose
pullback to $\mathcal{G}$ restricts to an isomorphism of each tangent space $T_p\mathcal{G}$ with $\mathfrak{g}$: this
pullback is a~generalization of the Maurer--Cartan form of $G$, and $M$ is locally isomorphic to~$G/P$ if\/f $\omega$ is
a~f\/lat connection.

R.~Baston and M.~Eastwood~\cite{BaEa:pt} have developed a~general theory of the Penrose transform for correspondences
$G/R\xleftarrow{\eta} G/Q\xrightarrow{\tau} G/P$ between generalized f\/lag varieties, where \mbox{$Q=P\cap R$}, and their theory
makes sense in the more general context of parabolic geometries studied by A.~\v Cap~\cite{Cap:cts}.
In this context, $G/P$ is replaced by a~parabolic geometry $N=\mathcal{G}/P$.
\v Cap shows that $M=\mathcal{G}/Q$ is then also a~parabolic geometry, and characterizes when it f\/ibres over
a~generalized twistor space $Y$ such that $\mathcal{G}$ is a~local principal $R$-bundle over~$Y$.
The Baston--Eastwood theory can be applied here, because the main ingredients of their construction only use the local
homogeneity of the f\/ibres of $\eta$ and $\tau$: in \v Cap's setting, the f\/ibres of $\eta\colon M\to Y$ are isomorphic to
open subsets of $R/Q$, which may be assumed contractible, while the f\/ibres of $\tau\colon M\to N$ are naturally
isomorphic to $P/Q$.
For any locally homogeneous vector bundle $\mathcal{E}$ on $Y$ (i.e., locally of the form $\mathcal{G}\times_R
V$, where $V$ is a~representation of $R$), the theory of~\cite{BaEa:pt} then provides a~spectral sequence
$E^{p,q}_1=\Gamma(N,\tau^q_*\Delta^p_\eta)$ which (assuming $N$ is Stein) converges to $H^{p+q}(Y,\mathcal{E})$: here
$\Delta^p_\eta$ is the relative (i.e., f\/ibrewise, for $R/Q$) Bernstein--Gelfand--Gelfand (BGG) resolution of
$\eta^{-1}\mathcal{E}$, and the direct image $\tau^q_*$ of this resolution can be computed using the Bott--Borel--Weil
theorem for $P/Q$.

This applies straightforwardly to projective surfaces and their minitwistor spaces, with $G=\SL(3,{\mathbb{C}})$, $P$
and $R$ the stabilizers of a~point and a~line in~${\mathbb{C}}{\rm P}^2$ respectively, where the point lies on the
line, so that $Q=P\cap R$ is the stabilizer of a~contact element.
A projective surface~$N$ has a~natural normal Cartan connection~$(\mathcal{G},\omega)$ and $M=\mathcal{G}/Q$ is the
correspondence space~$P(TN)$.
\v Cap's criterion for the local $R$-action of~$\mathcal{G}$ over $Y$ is automatically satisf\/ied because the projective
spray has rank one.
The relevant locally homogeneous bundles on~$Y$ are the bundles $\mathcal{O}_Y(k) S^\ell TY$.
The computation of the Penrose transforms are straightforward following~\cite{BaEa:pt}, since the generalized f\/lag
varieties $P/Q$ and $R/Q$ are both projective lines, and the spectral sequence collapses at $E^{p,q}_2$.

{\sloppy The results can be interpreted in terms of BGG sequences of projectively invariant li\-near dif\/ferential operators
(PILDOs).
These are `curved versions' of the BGG resolutions for $\SL(3,{\mathbb{C}})/P$ (in this case) which were shown to exist
(in general) by \v Cap, Slov\'ak and Sou\v cek~\cite{CSS:bgg}.
In the projective case, the BGG sequences consist of two PILDOs associated to each irreducible representation of
$\SL(3,{\mathbb{C}})$.
I~ask the reader's indulgence for not writing out all the calculations, and give the results of the analysis without
proof.

}

\begin{prop}
Assuming $N$ is Stein and geodesically convex, the Penrose transforms of the cohomology of
$\mathcal{E}=\mathcal{O}_Y(k-\ell)\otimes S^\ell TY$ for $\ell\in\{0,1\}$ are given as follows, where $\odot$ denotes
Cartan product of representations of $\SL(3,{\mathbb{C}})$ $($which amounts to taking the tracefree part of the
tensor product in the examples below$)$.
\begin{itemize}\itemsep=0pt
\item
If $k\geqslant 0$ then $H^2(Y,\mathcal{E})=0$, whereas $H^0(Y,\mathcal{E})$ and $H^1(Y,\mathcal{E})$ are the kernel and
cokernel of the first BGG operator for $S^\ell{\mathbb{C}}^{3*}\odot S^{k}{\mathbb{C}}^{3}$, which is a~PILDO of order
$\ell+1$
\begin{gather*}
\mathcal{O}_N(2k+\ell) S^{k}T^*N\to\mathcal{O}_N(2k+\ell) S^{k+\ell+1}T^*N.
\end{gather*}
\item
If $k\leqslant -1$ and $k+\ell\geqslant -2$ then $H^0(Y,\mathcal{E})=0$ and $H^2(Y,\mathcal{E})=0$, whereas
$H^1(Y,\mathcal{E})$ is the space of sections of the bundle
\begin{gather*}
\mathcal{O}_N(2k+\ell+3) S^{-k-2}TN\oplus\mathcal{O}_N(2k+\ell)S^{k+\ell+1}T^*N
\end{gather*}
$($with the convention that $S^{-1}TN=0)$.
\item
If $k+\ell\leqslant -3$ then $H^0(Y,\mathcal{E})=0$, whereas $H^1(Y,\mathcal{E})$ and $H^2(Y,\mathcal{E})$ are the
kernel and cokernel of the second BGG operator for $S^{-(k+\ell)-3}{\mathbb{C}}^{3*}\odot S^{\ell}{\mathbb{C}}^{3}$,
which is a~PILDO of order $\ell+1$
\begin{gather*}
\mathcal{O}_N(2k+\ell+3)  S^{-k-2}TN\to\mathcal{O}_N(2k+\ell+3)  S^{-(k+\ell)-3}TN.
\end{gather*}
\end{itemize}
For $\ell\geqslant 2$ the same results hold, except that the identification of the differential operators with BGG
operators may hold only modulo lower order terms involving the curvature $C$ of the projective structure.
\end{prop}
(It is possible that a~more concrete analysis could identify the PILDOs with the BGG operators in all cases, but this
seems to require understanding the normal projective Cartan connection in twistorial terms, which does not appear
entirely straightforward to me.)

Some of these Penrose transforms are easy to see by other means.
For example, when \mbox{$k=\ell=1$}, $H^0(Y,TY)$ is obviously the space of projective vector f\/ields, which is the kernel of
a~second order dif\/ferential operator on sections of $TN\cong\mathcal{O}_N(3)T^{*} N$, and the f\/irst BGG operator for the
adjoint representation ${\mathbb{C}}^{3*}\odot{\mathbb{C}}^3$.

The case $k=3,\ell=0$ is interesting from the point of view of~\cite{Hit:gse}: in this case the bundle over $Y$ is the
anticanonical bundle and its sections are Poisson structures.
These correspond, via Penrose transform, to sections of $\mathcal{O}_N(6) S^3T^{*} N\cong S^2TN\otimes_0 T^{*} N$ in the
kernel of the natural f\/irst order PILDO with values in $\mathcal{O}_N(6) S^4T^{*} N\cong S^2TN\otimes_0 S^2T^{*} N$.

\subsection*{Ward correspondences}

The Ward correspondence (cf.~\cite{Ward:sdg}) for degree zero bundles over the minitwistor space of a~projective
surfaces is completely straightforward.
It applies to vector bundles or principal bundles (with af\/f\/ine algebraic structure group) over $Y$ which are trivial on
minitwistor lines, or, more generally, to f\/ibre bundles $\mathcal{E}\to Y$ such that $\mathcal{E}|_{\hat x} \cong
\hat x\times H^0(\hat x,\mathcal{E})$ for any minitwistor line $\hat x$.
Let~$\mathcal{E}$ be such a~f\/ibre bundle and let $E$ be the f\/ibre bundle over $N$ whose f\/ibre at $x\in N$ is $H^0(\hat
x,\mathcal{E})$.
This is equipped with a~connection whose parallel sections along geodesics are given by elements of the corresponding
f\/ibre of~$\mathcal{E}$.
Conversely, any connection on a~f\/ibre bundle $E$ over $N$ is f\/lat along any geodesic, and hence the spaces of parallel
sections along geodesics def\/ine a~f\/ibre bundle $\mathcal{E}$ over $Y$.
This f\/ibre bundle is trivialized along $\hat x$ by the parallel extensions of any nonzero element of $E_x$, and this
identif\/ies $\mathcal{E}|_{\hat x}$ with $\hat x\times H^0(\hat x,\mathcal{E})$ provided that any holomorphic map from
${\mathbb{C}}{\rm P}^1$ to $E_x$ is constant, which certainly holds if $E_x$ is a~vector space or an af\/f\/ine
algebraic variety.

\begin{prop}
There is a~bijective correspondence between vector bundles with connection on $N$ and vector bundles on $Y$ which are
trivial on minitwistor lines.
\end{prop}

Notice that, in the case of (degree zero) line bundles, this construction reduces to the Penrose transform of
$H^1(Y,\mathcal{O}_Y)$, which is the cokernel of ${\rm d}$ from functions to $1$-forms.

Given a~vector bundle $\mathcal{E}$ which is trivial on minitwistor lines, the Penrose transform above may be
generalized to cohomology with values in the tensor product of a~locally homogeneous bundle with $\mathcal{E}$ (or any
bundle constructed from $\mathcal{E}$): one simply couples the corresponding dif\/ferential operators to the connection on
$E$ (these operators are all strongly invariant~\cite{CSS:bgg} and so can be coupled naturally to connections).
The same idea applies to principal bundles via the associated bundle construction.

This provides in particular, a~Ward correspondence for projective pairs $(\nabla^\alpha,\phi)$.
I will assume the gauge group acts faithfully and transitively on an af\/f\/ine algebraic variety $\Sigma$.
Then the bundle $E\to N$ associated to this action has a~connection $\nabla^\alpha$, which the Ward correspondence
associates with a~bundle $\mathcal{E}\to Y$ trivial on minitwistor lines.
By Penrose transform, $\phi$ is in the kernel of the PILDO $\sym D$ on $\mathcal{O}_N(2) T^{*} N$ coupled to $\nabla$ on
the gauge algebra of vector f\/ields.
Hence it corresponds, via Penrose transform, to a~vertical vector f\/ield~$V$ on~$\mathcal{E}$ with values in
$\mathcal{O}_Y(1)$.
If~$V$ is nonvanishing, then there is a~quotient bundle $\mathcal{F}$ over~$Y$ which is not trivial on minitwistor
lines: instead its vertical bundle has degree one.

\subsection*{Divisors and congruences}

The third ingredient of any respectable twistor correspondence is an analysis of the geometric objects corresponding to
divisors in the twistor space.
In this case a~degree one divisor $\mathcal{C}$ in $Y$ (i.e., $\mathcal{C}$ meets each minitwistor line in
a~point, so that its divisor line bundle has degree one) corresponds to a~geodesic congruence in $N$.

The divisor line bundle of $\mathcal{C}$ may be written $\mathcal{L}_\mathcal{C}\otimes \mathcal{O}_Y(1)$ where
$\mathcal{L}_\mathcal{C}$ has degree zero, hence corresponds to a~line bundle $L_\mathcal{C}$ on $N$ with connection
$\nabla$.
The (unique up to scale) nonzero section~$s$ of $\mathcal{L}_\mathcal{C}\otimes \mathcal{O}_Y(1)$ vanishing on
$\mathcal{C}$ corresponds, via Penrose transform, to a~nonzero section~$\phi$ of $L_\mathcal{C}\otimes \mathcal{O}_N(-1)
TN$ satisfying $D^\nabla\phi=\frac12 {\rm d}^{D,\nabla}\phi$ for any connection $D$ in the projective class, and
$\phi$ is tangent to the geodesic congruence because $s$ vanishes on~$\mathcal{C}$; notice also that $\nabla$ is the
unique connection on $L_\mathcal{C}$ with the property that $D^\nabla\phi$ is skew, and it is f\/lat if\/f $\mathcal{C}$ is
(locally) a~divisor for $\mathcal{O}_Y(1)$.
Also, if $D$ is a~connection in the projective class with $D^\nabla\phi=0$, then
$\abrack{{\rm r}^D\wedge{\rm id}}\cdot\phi+F^\nabla\otimes\phi=0$; in particular,
${\rm r}^D(\phi,\phi)=0$.

The theory of degree two divisors is more interesting.
Locally and generically, such a~divisor is the sum $\mathcal{C}=\mathcal{C}_1+\mathcal{C}_2$ of two degree one divisors
which correspond to a~pair of linearly independent geodesic congruences $\phi_j$ ($j=1,2$), these being sections of
$L_{\mathcal{C}_j}\otimes \mathcal{O}_N(-1) TN\cong L_{\mathcal{C}_j}\otimes \mathcal{O}_N(2) T^{*} N$ for line bundles
$L_{\mathcal{C}_j}$ with connections $\nabla^j$.
The images of $\phi_1$ and $\phi_2$ are null directions for the conformal metric $\mathsf{c}=\phi_1\odot\phi_2$, which
is the section of $L_{\mathcal{C}_1}\otimes L_{\mathcal{C}_1}\otimes\mathcal{O}_N(4) S^2T^{*} N$ obtained from the
$\mathcal{L}_\mathcal{C}$-coupled Penrose transform, where
$\mathcal{L}_{\mathcal{C}}=\mathcal{L}_{\mathcal{C}_1}\otimes\mathcal{L}_{\mathcal{C}_2}$.
Thus $\sym D^\nabla\mathsf{c}=0$, where $\nabla$ is the tensor product connection on $\mathcal{L}_{\mathcal{C}}$.

The area form (determinant) of $\mathsf{c}$ identif\/ies $L_{\mathcal{C}_1}\otimes L_{\mathcal{C}_2}$ with
$\mathcal{O}_N(-1)$, so that $\mathsf{c}$ is a~genuine metric on $\mathcal{O}_N(-3/2) TN$.
There is now a~unique connection $D$ in the projective class with $D\mathsf{c}=0$, so that $D$ is the Weyl connection of
$\mathsf{c}$.
Since $D$ preserves the null directions, it follows easily that $D^{\nabla^j}\phi_j=0$ and on $\mathcal{O}_N(-1)$, and
by the characterization of $\nabla^1$ and $\nabla^2$, $D=\nabla$ on $\mathcal{O}_N(-1)$.
Thus $\abrack{{\rm r}^D\wedge{\rm id}}\cdot\phi_j+F^{\nabla^j}\otimes\phi_j=0$ for $j=1,2$, and
straightforward computation yields the following result, which (as usual) should be interpreted locally (to
a~neighbourhood of a~point in $N$ or a~minitwistor line in $Y$).

\begin{prop}
Let $\mathcal{C}=\mathcal{C}_1+\mathcal{C}_2$ be a~degree two divisor on $Y$ corresponding to a~conformal metric
$\mathsf{c}$ with Weyl connection $D$, and geodesic congruences $\phi_1$ and $\phi_2$ with associated connections
$\nabla^1$ and $\nabla^2$.
\begin{itemize}\itemsep=0pt
\item
The following are equivalent:
\begin{enumerate}\itemsep=0pt
\item[$(i)$]
$D$ is flat on $\mathcal{O}_N(1)$;
\item[$(ii)$]
${\rm r}^D$ is symmetric;
\item[$(iii)$]
$F^{\nabla^1}+F^{\nabla^2}=0$;
\item[$(iv)$]
$\mathcal{C}$ is a~divisor for $\mathcal{O}_Y(2)$.
\end{enumerate}
The $D$-parallel sections of $\mathcal{O}_N(1)$ determine, with $\mathsf{c}$, a~homothetic family of metrics with
Levi-Civita connection $D$ compatible with the projective structure.
\item
The following are also equivalent:
\begin{enumerate}\itemsep=0pt
\item[$(i)$]
$D$ is flat on $\mathcal{O}_N(-3/2) TN$;
\item[$(ii)$]
${\rm r}^D$ is skew;
\item[$(iii)$]
$F^{\nabla^1}-F^{\nabla^2}=0$;
\item[$(iv)$]
$\mathcal{C}_1-\mathcal{C}_2$ is the divisor of a~meromorphic function, hence $Y$ fibres over
${\mathbb{C}}{\rm P}^1$.
\end{enumerate}
The $D$-parallel sections of $\mathcal{O}_N(-3/2) TN$ then correspond to the fibres of $Y\to{\mathbb{C}}{\rm P}^1$.
\end{itemize}

Combining the two equivalences, it follows that $Y$ admits distinct divisors $\mathcal{C}_1$ and $\mathcal{C}_2$ for
$\mathcal{O}_Y(1)$ if and only if the projective structure is flat.
\end{prop}

The case of $Y$ f\/ibering over ${\mathbb{C}}{\rm P}^1$ is considered by Dunajski and West~\cite{DuWe:acs}.
This condition can be interpreted on the correspondence space as the existence of a~trivialization with respect to which
the projective spray has no $\partial_\lambda$ term, and they show that this means that the ODE associated to the
projective structure is point equivalent to the derivative of a~f\/irst order ODE.
On the other hand, the case that $Y$ has a~canonical degree two divisor $\mathcal{C}$ arises in the
characterization~\cite{LeMa:zmcs}, by LeBrun and Mason, of Zoll metrics among Zoll projective structures: in their
setting, $Y={\mathbb{C}}{\rm P}^2$ and so such a~degree two divisor $\mathcal{C}$ is a~conic.

Degree three divisors include the case of sections of the anticanonical bundle mentioned already in this appendix, but I
have nothing more to say about this.

\subsection*{Acknowledgements}

I am extremely grateful to Maciej Dunajski and Simon West for introducing me to their sti\-mu\-lating work, and for several
helpful comments. I also thank the EPSRC for f\/inancial support in the form of an Advanced Research Fellowship.

\pdfbookmark[1]{References}{ref}
\LastPageEnding

\end{document}